\documentclass[12pt]{article}
\usepackage{amsfonts, amssymb, amsthm}

\newcommand{\sect}[1]{\setcounter{equation}{0}\section{#1}}
\newcommand{\subsect}[1]{\subsection{#1}}

\def\A{{\cal A}}
\def\B{{\cal B}}
\def\C{{\cal C}}
\def\z{{\hat z}}

\def\co{{\Delta}}

\usepackage{times}

\begin{document}
%%%%%%%%%%%%%%%%%%%%%%%%%%%%%%%%%%%%%%
\thispagestyle{empty}
\hfill\today
\vspace{2.5cm}

\begin{center}
{\LARGE{\bf{Three dimensional quantum algebras:
\\[0.45cm]
a Cartan-like point of view}}}
\end{center}

\bigskip\bigskip

\begin{center}
A. Ballesteros$^1$, E. Celeghini$^2$  and M.A. del Olmo$^3$
\end{center}

\begin{center}
$^1${\sl Departamento de F\'{\i}sica, Universidad de Burgos, \\
E-09006, Burgos, Spain.}\\
\medskip

$^2${\sl Departimento di Fisica, Universit\'a  di Firenze and
INFN--Sezione di
Firenze \\
I50019 Sesto Fiorentino,  Firenze, Italy}\\
\medskip

$^3${\sl Departamento de F\'{\i}sica Te\'orica, Universidad de
Valladolid, \\
E-47011, Valladolid, Spain.}\\
\medskip

{e-mail: angelb@ubu.es, celeghini@fi.infn.it, olmo@fta.uva.es}
\end{center}

\bigskip

\begin{abstract}
A perturbative quantization procedure for Lie bialgebras is introduced
and used to classify all three
dimensional complex quantum algebras compatible with a given coproduct.
The role of elements of the
quantum universal enveloping algebra that, analogously to generators in
Lie algebras, have a
distinguished type of coproduct is discussed, and the relevance of a
symmetrical basis in the
universal enveloping algebra  stressed. New 
quantizations
of three dimensional solvable  algebras, relevant for possible physical applications for their
simplicity, are obtained  and all already known related results
recovered. Our results give a
quantization of all existing three dimensional Lie algebras and
reproduce, in the classical limit, the
most relevant sector of the complete
classification for real three dimensional Lie bialgebra structures given
in
\cite{gomez}.
\end{abstract}
\vskip 1cm

\vskip 0.4cm

\vfill
\eject
%%%%%%%%%%%%%%%%%%%%% INTRODUCTION %%%%%%%%%%%%%%%%%%%%%%%%%%%%%%%%%
\sect{Introduction\label{introduccion}}

Cartan classification of semisimple Lie algebras has facilitated 
their applications in Physics. Quantum
algebras are not been classified in a similar way and their physical applications are far to be
developed or understanding. For these reasons, 
this paper deals with some facets of the problem of the construction and
classification
of quantum universal enveloping algebras (hereafter,
quantum algebras) \cite{Dri}-\cite{Tjin}. It is well-known that any
quantum algebra
$(U_z(\mathfrak a),\Delta)$ with deformation parameter $z$ defines a
unique Lie bialgebra structure
$(\mathfrak a,\eta)$, a pair determined by the Lie
algebra
$\mathfrak a$  and a  skew-symmetric linear map  (cocommutator)
$\eta:\mathfrak a \to \mathfrak a\times
\mathfrak a$ . Such cocommutator $\eta$ is defined as the first
order skewsymmetric part of the coproduct
\begin{equation}\label{bialgebra}
\eta (X)=\frac 1 2 (\Delta (X)-\sigma\circ \Delta (X))\quad
\mbox{mod}\,z^2\qquad\forall X\in
\mathfrak a
\end{equation}
where $ \sigma$ is the flip operator  $ \sigma (X\otimes
Y)=Y\otimes X$.

Therefore, quantum deformations of a
given Lie algebra
$\mathfrak a$ can be classified according to this
``semiclassical" limit. Moreover, we recall that Lie bialgebras
are in one to one correspondence with Poisson-Lie structures on the
group
$Lie(\mathfrak a)$ \cite{DriPL} that arise again as the first order in
$z$ of the
quantum groups dual to the quantum algebras
$(U_z(\mathfrak a),\Delta)$.
With this in mind, some classifications of Lie bialgebra structures for
several
physically relevant Lie algebras have been obtained. For the three
dimensional (3D) case  we will refer to the full classification of Lie
bialgebras given in \cite{gomez} and for some higher dimensional Lie
bialgebras we refer to
\cite{galileo} and references therein.

However, it is clear that the inverse problem has also to be
 faced: that is, to obtain general recipes for the ``quantization" (i.e.
the associated quantum
algebra) of a given Lie bialgebra
$(\mathfrak a,\eta)$. Although the existence of such an object is indeed
guaranteed (see
\cite{charipressley}, Chapter 6), only for coboundary triangular
bialgebras the Drinfel'd twist
operator gives rise to the associated quantum algebra
\cite{majid}. However, quasitriangular and even non-coboundary Lie
bialgebras do exist and for them
the twist operator approach to quantization is not available. We recall
that some attempts have been
performed in order to get structural properties of the quantization of
arbitrary Lie bialgebras (see
\cite{LM,mudrov} for a prescription to get the quantum coproduct --but
not the deformed commutation
rules-- for a wide class of examples). Moreover, to our knowledge, a
general investigation concerning
the uniqueness of this Lie bialgebra quantization process has not been
given yet and only
restrictive results for certain deformations of simple Lie algebras have
been obtained (see
\cite{schnider}, Chapter 11). As a consequence of the abovementioned
facts, a complete
classification of quantum algebras in the spirit of the Cartan
classification for Lie algebras is still far to be reached.

In this paper we present  a
``direct approach" to the quantization problem, together with an
operational procedure {\it \`a la
Cartan} (in a sense that will be made more explicit in the sequel) for
the classification of quantum
algebras. The essential ingredients of the quantization method are
described in Section 2. In
particular,  given a Lie bialgebra we shall firstly obtain a
coassociative quantum coproduct with some
outstanding symmetry properties and, afterwards, we shall solve order by
order the compatibility
equations for the deformed commutation rules. Lie-Poisson brackets 
are immediately recover 
using a
completely symmetrized basis for
$U_z(\mathfrak a)$  non  previously considered in the
literature. In Section 3 this
approach will be explicitly developed for a relevant type of
cocommutator $\eta$ that is
compatible with all the non-isomophic 3D complex Lie
algebras, as described -for
instance- in Jacobson \cite{Jac}. This cocommutator 
gives rise to a coproduct very simple and 
deformed commutators that are not specially complex.
In this way we obtain three different
families of quantum algebras
whose deformed commutation rules are very general and depend on several
structure constants.

Another point to stress is that our quantization method has been made using a symmetrized basis
instead the usual Poincar\'e-Birkhoff-Witt (PBW)  basis. 
This procedure has not been previously
considered in the literature and
can be relevant in quantum physics
applications and in the semiclassical limit.

In order to classify such quantizations, in Section 4 we consider the
equivalence
transformations on the basis of the quantum algebra that are defined
through invertible maps in
$U_z(\mathfrak a)$ that leave formally invariant the coproduct for the
elements of the basis.  This
constraint regarding the formal invariance of the  coproduct can be
understood as a way to identify
certain basis elements as some sort of ``generators" of a quantum
algebra, thus following for quantum
algebras the ``Cartan" approach to the classification of Lie algebras
within
$U(\mathfrak a)$, where only linear transformations in the space of
generators leaving invariant their
(primitive) coproduct are performed. After a such classification is
systematically developed, new
quantizations of 3D solvable algebras are obtained, and  already
known results are recovered. The above mentioned simplicity of the 
algebraic and coalgebraic structures of this new quantizations increases notably their possible 
interest
in Physics . In order to make more precise the range of 3D deformations that have been
covered, a detailed comparison
with the complete classification of 3D real Lie bialgebra
structures presented in
\cite{gomez} is explicitly given within a final Section which also
includes some further comments and
conclusions.

%%%%%%%%%%%%%%%%%%%%%%%%%%%%%%%%%%%%%%%%%%%%%%%%%%%%%%%%%%%%%%%
%%%%%%%%%%%%%%%%%%%%%%%%%%%%%%%%%%%%%%%%%%%%%%%%%%%%%%%%%%%%%%%
\sect{The quantization method}\label{generalquantization}

Let us consider the Lie algebra
$\mathfrak a$ and its universal enveloping algebra $U(\mathfrak a)$
\cite{Jac}, an associative algebra that is obtained as the quotient
$T(\mathfrak a)/I$, where $T(\mathfrak a)$ is the tensor algebra of
$\mathfrak a$ and
$I$ is the ideal generated by the elements $XY-YX-[X,Y]$ ($X,Y\,\in
\mathfrak a$). If we define the coproduct, counit and antipode
($\forall\, X \,\in  \mathfrak a$)
$$
\co_0(X)=1\otimes X + X\otimes 1,\quad \co_0(1)=1\otimes 1,\quad
\epsilon(X)=0,\quad \epsilon(1)=1,\quad
\gamma(X)=-X,
$$
and we extend by linearity all these (anti)automorphisms to the full 
$U(\mathfrak a)$, we shall endow the universal enveloping algebra with a
Hopf algebra structure. In general, an element
$Y$ of a Hopf algebra is called {\it primitive} if $$\co(Y)=1\otimes Y +
Y\otimes 1.$$ 
Within $U(\mathfrak a)$, it can be shown that the only
primitive generators
under the coproduct  $\co_0$ are the generators of
$\mathfrak a$ (this result is known as the Fiedrichs theorem
\cite{Post}). Note that an essential property of the Hopf algebra
$U(\mathfrak a)$ is its cocommutativity, since $\co_0$ is invariant
under the action
of the flip operator $\sigma$.

Let us now consider the Hopf algebra $U(
\mathfrak a)$ and a deformation parameter $z$.
A quantum algebra $(U_z(\mathfrak a),\co)$  is a Hopf algebra of formal
power series in the deformation parameter $z$ with coefficients in $U(
\mathfrak a)$ and such that \cite{charipressley}
$$
U_z(\mathfrak a)/zU_z(\mathfrak a)\simeq U(\mathfrak a).
$$
Therefore, $U(\mathfrak a)$ is obtained (as Hopf algebra) in the
limit $z\to 0$, and the first order in $z$ of the
coproduct is directly related to the cocommutator of an underlying Lie
bialgebra $(\mathfrak a,\eta)$ through (\ref{bialgebra}). In this way,
Lie bialgebras can be
used to characterize quantum deformations. Amongst all the Hopf algebra
axioms to be imposed,
we recall that the quantum coproduct
$\Delta$ has to be a coassociative map, namely:
\begin{equation}\label{coassociative}
(\Delta\otimes id)\circ\Delta=(id\otimes \Delta)\circ\Delta.
\end{equation}

From now on, we shall refer to the ``quantization" of a given Lie
bialgebra $(\mathfrak a,\eta)$ as the
problem of finding a quantum algebra $(U_z(\mathfrak a),\co)$ such that
(\ref{bialgebra}) is fulfilled
(we recall that the uniqueness of such construction cannot be take for
granted).

The general quantization procedure that we propose is based on three
essential ingredients that can be
considered whatever the dimension of the Lie bialgebra is. The first one
is a generalized
cocommutativity property for the coproduct, the second one is related to
the choice of a basis 
in the universal enveloping algebra and the third one selects a given
type of power series expansion for the deformed commutation rules. 
%%%%%%%%%%%%%%%%%%%%%%%%%%%%%%%%
%%%%%%%%%%%%%%%%%%%%%%%%%%%%%%%% 

\medskip
\noindent {\bf 1. Generalized cocommutativity}. In general, we impose
the (noncocommutative)
quantum coproduct $\Delta$
to be invariant under the composition $\tilde \sigma=\sigma\circ T$ of
the flip operator $ \sigma $
and a change  of sign of (all) the deformation parameter(s):$$
\tilde \sigma
\circ \Delta  = \Delta \qquad \mbox{where}\quad  \tilde
\sigma=T\circ \sigma \quad \mbox{and}\quad T(z)=-z. $$
We point out that the definition (\ref{bialgebra}) for the underlying
cocommutator implies that the
deformation parameter appear explicitly as multiplicative factors
within the cocommutator. This
symmetry property of the coproduct can be imposed in any dimension and
makes much easier the
procedures of symmetrization and `hermitation'
\cite{enrico03}. In particular, given a certain Lie bialgebra
$(\mathfrak a,\eta)$,
the above assumption implies that the first order deformation of the
coproduct will be just given by the
(skewsymmetric) Lie bialgebra cocommutator.
\begin{equation}\label{first}
\Delta(X)=\Delta_0 (X) + \eta(X) + O[z^2], \qquad X\in\mathfrak{g}.
\end{equation}

Moreover, the invariance of $\Delta$
under the transformation $\tilde \sigma$ together with the fact that the
coproduct is an
algebra homomorphism with respect to the deformed commutation rules
$[\cdot,\cdot]_z$, i. e.
$\Delta ([\cdot,\cdot]_z)=[\Delta (\cdot), \Delta (\cdot)]_z$ implies
that $[\cdot,\cdot]_z$ has to
be an even function in the deformation parameter $z$.

%%%%%%%%%%%%%%%%%%%%%%%%%%%%
%%%%%%%%%%%%%%%%%%%%%%%%%%%%

\medskip
\noindent {\bf 2. The choice of a basis in $U_z(\mathfrak a)$}. In
contradistinction with previous
works on this subject in which the PBW   basis
$X_1^\alpha\,X_2^\beta\,\dots\,X_l^\zeta$ is considered,
we introduce using the operator $\rm Sym$ a basis in
$U_z(\mathfrak a)$ given by the completely  symmetrized monomials.
We, thus, define the linear operator ${\rm Sym}$  by
$$
{\rm Sym}\;\left\{\sum c_{\alpha\beta\dots\zeta}\;
X_1^\alpha\,X_2^\beta\,\dots\,X_l^\zeta\right\} =
\sum c_{\alpha\beta\dots\zeta}\, {\rm
Sym}\,\{X_1^\alpha\,X_2^\beta\,\dots\,X_l^\zeta\}
$$ 
$$
{\rm Sym}\;\{A_1\dots A_n\} := \frac 1{n!}\sum_{p\in {\rm S_n}} p (A_1\dots A_n),
$$
with ${\rm S_n}$  the group of permutations of n elements .
Note that Sym  is the identity for commuting operators.

 This
symmetrization Ansatz, although very convenient in quantum mechanical
terms, has not been previously
considered in the literature, and turns out to be very efficient in
order to get the explicit form
of the deformed commutation rules.
We remark that one of the main advantages of the symmetrized basis in
the quantum algebra is the fact
that if we replace the deformed commutation rules by Poisson brackets,
the correspondent Poisson-Hopf
algebra is uniquely and immediately defined.

%%%%%%%%%%%%%%%%%%%%%%%%%%%%
%%%%%%%%%%%%%%%%%%%%%%%%%%%%

\medskip
\noindent {\bf 3. Deformed commutation rules}. We will assume
\begin{equation}\label{hipotesis3}
[X,Y]_z = \frac{1}{z} {\rm Sym}\;(f( 
zX_1,zX_2,\dots ,zX_l)),   
\end{equation} 
where $f$ is a meromorphic function at $z=0$ and also odd in $z$; 

We extend the discussion to 
the multiparametric case considering  meromorphic deformations 
of Lie algebras with $z_i/z_j$ fixed.

%Since all the quantum algebras look to us described by meromorphic functions 
%in the parameter $z$ around $z=0$ we extend the discussion to the multiparametric 
%case considering that
%$z_i/z_j$ fixed.

 Under all these assumptions, given a (family of) Lie bialgebra
$(\mathfrak a,\eta)$, the ``direct"
quantization procedure that we propose would be sketched as follows:
\begin{itemize}

\item Assume that the first order coproduct is of the form
(\ref{first}).

\item Order by order in the deformation parameter(s), get the relations
coming from the
coassociativity constraint (\ref{coassociative})
and solve them recursively by taking into account the
invariance under
$\tilde\sigma$ of the solution, thus obtaining the full quantum
coproduct.

\item Obtain, again order by order, the deformed commutation rules by
solving the compatibility
equations coming from the fact that the coproduct has to be an algebra
homomorphism.

\end{itemize}

%%%%%%%%%%%%%%%%%%%%%%%%%%%%%%%%%%%%%%%%%%%%%%%%%%
%%%%%%%%%%%%%%%%%%%%%%%%%%%%%%%%%%%%%%%%%%%%%%%%%%%%%%%%
\sect{Three dimensional quantum algebras}

Let $\{ A,B,C\}$ be the generators of an arbitrary complex
3D  Lie
algebra ${\cal L}$ with commutators
\begin{equation}\label{conmutaciones}
\begin{array}{lll}
[A,B] &=& c_1 A + c_2 B + c_3 C  \\[0.3cm]
[A,C] &=& b_1 A + b_2 B + b_3 C \\[0.3cm]
[B,C] &=& a_1 A + a_2 B + a_3 C ,
\end{array}\end{equation}
where the structure constants are complex numbers subjected to some
(nonlinear) relations coming from
 Jacobi identity.

We recall that the complete classification  of the 3D
complex Lie algebras is given in
\cite{Jac} (see for instance \cite{pavel} for the real case).
According to the dimension of the
derived  algebra
${\cal L}'=[{\cal L},{\cal L}]$, the
 nonisomorphic classes of 3D Lie algebras read:
\begin{itemize}
\item Type I: ${\rm dim}\;{\cal L}'=0$. Then ${\cal L}$ is abelian. We
shall not consider this case
from the point of view of quantum deformations, since any coassociative
coalgebra with primitive
non-deformed limit is compatible with the abelian commutation rules.

\item  Type II:  ${\rm dim}\;{\cal L}'=1$.  We have two algebras,  the
Heisenberg-Weyl algebra  and a
central extension ${\cal L}={\cal B}\oplus{\cal C}$ of the Borel
algebra, where $ {\cal B}$ is a Borel
algebra and ${\cal C}$ commutes with
${\cal B}$.

\item Type III:  ${\rm dim}\;{\cal L}'=2$.  We have the family on
non-isomorphic Lie algebras labeled
by the nonzero complex number $\alpha$ with commutators
\begin{equation}\label{jacobson1}
[X_1,X_2]=0,\qquad [X_1,X_3]=X_1,\qquad
[X_2,X_3]=\alpha X_2 ,
\end{equation}
and the Lie algebra
\begin{equation}\label{32}
[X_1,X_2]=0,\qquad [X_1,X_3]=X_1+X_2,\qquad
[X_2,X_3]= X_2 .
\end{equation}

\item  Type IV: ${\rm dim}\;{\cal L}'=3$. The only element in this class
is the simple Lie algebra
${\cal A}_1$.
\end{itemize}

Note that in the Type III   the algebra of dilations and the 2D Euclidean
algebra for $\alpha=1$ and $\alpha=-1$, respectively, are included.

 A complete classification of the Lie bialgebra structures of
the {\it real} 3D Lie
algebras has been given in \cite{gomez}, and the Lie bialgebra
classification for the complex case can
be extracted from there.

Obviously, once a set of values for the structure constants
$\{a_i,b_i,c_i\}$  is given, a suitable linear
transformation $X_i=X_i(A,B,C)$ with
complex coefficients can be found in such a way that the Lie algebra
(\ref{conmutaciones}) is reduced
to one of the Jacobson cases.
 However, while in (\ref{conmutaciones}) all $X_i$ are primitive and, 
thus, equivalent, this is not longer
true in presence of deformation. 

 In particular, by looking at 
the complexified form of the
classification of 3D Lie
bialgebras  given in
\cite{gomez}, one can realize that a cocommutator of the form
\begin{equation}\label{coconmutador}
\eta(A)=0  \qquad
\eta(B)= z\, A \wedge B ,  \qquad
\eta(C)=\rho\,z\,  A \wedge  C , \qquad
z,\rho\in \mathbb C
\end{equation}
defines a Lie bialgebra structure for each of the types of Lie algebras
given in Jacobson's
classification provided that different linear transformations
$X_i=X_i(A,B,C)$ are defined and that
$\rho$ takes certain appropriate values. Since we are interested in
obtaining the most general types of deformed commutation rules arising
in 3D quantizations, we shall
apply the perturbative quantization procedure described in the previous
Section in order to get a
quantum coproduct coming from the cocommutator (\ref{coconmutador})
together with a
compatible deformation of the commutation rules (\ref{conmutaciones}).
We remark that the cocommutator (\ref{coconmutador}) can be thought of
as a two parameter
structure in
$z$ and  $\chi=\rho\,z$. Moreover, complex values for all the parameters
(including the
structure constants $\{a_i,b_i,c_i\}$ and $\rho$) will be considered,
and the results here presented
will also contain the quantizations in which the deformation parameter
is a root of unity.

Let us follow step by step the procedure introduced in Section 2.
Firstly, we assume that the quantum
coproduct will be of the form (\ref{first}), namely
\begin{equation}\label{assumption}
\begin{array}{lll}
\Delta(A) &=&  1\otimes  A +   A \otimes 1 , \\[0.3cm]
\Delta(B) &=&  1\otimes  B +   B \otimes 1 +  z\, A \wedge B + O[z^2],
\\[0.3cm]
\Delta(C) &=&  1\otimes  C +   C \otimes 1 +  \rho\,z\,  A \wedge  C  +
O[z^2].
\end{array}\end{equation}
In this case, it is straightforward to prove that the following
coassociative coproduct has a first
order given by the  cocommutators (\ref{coconmutador}) and is invariant
under the transformation
$\tilde\sigma$:
\begin{equation}\label{coproductos}
\begin{array}{lll}
\Delta(A) &=&  1\otimes  A +   A \otimes 1 , \\[0.3cm]
\Delta(B) &=&  e^{z A}\otimes  B +   B \otimes e^{-z A} , \\[0.3cm]
\Delta(C) &=&  e^{ \rho\, z A}\otimes  C +   C \otimes e^{- \rho\, z A},
\end{array}\end{equation}

Now we have to obtain the deformed commutation rules by solving the
compatibility equations coming from the fact that the coproduct has to
be an algebra homomorphism.
Since the $\tilde\sigma$ invariance of the coproduct implies that the
first deformed term in the
commutation rules has to be of the order $z^2$, the first order
coproduct (\ref{assumption}) has to be
compatible with the {\it non deformed} commutation rules.
This condition leads to the following relations between the structure
constants and $\rho$:
\begin{equation}\label{constantesp}
b_2 (1-\rho )=0,\quad c_3 (1-\rho)=0, \quad
b_1=-\rho\,a_2,\quad   a_3=\rho\,c_1,\quad   b_3  = - \rho\, c_2.
\end{equation}
Jacobi identity requires the relation
\begin{equation}\label{jacobi0}
(1-\rho) [a_2 c_1(1+\rho)  -a_1 c_2] =0
\end{equation}
Now we can distinguish the following cases:

\begin{enumerate}

\item
Case $\rho \neq \pm 1$:
\begin{equation}
\begin{array}{c}\label{constantesp0}
b_2 =0,\qquad c_3=0, \quad
b_1=-\rho\,a_2,\quad   a_3=\rho\,c_1,\quad   b_3  = - \rho\, c_2    ;
\\[0.3cm]
a_2 c_1(1+\rho)  -a_1 c_2=0.
\end{array}\end{equation}
Note that the space of parameters is 3D
 
\item
Case $\rho =+1$:
\begin{equation}\label{constantesp1}
b_1=-\,a_2,\quad   a_3=\,c_1,\quad   b_3  = - \, c_2.
\end{equation}
In this case the space of parameters is 6D.
 
\item
Case $\rho =-1$:
\begin{equation}\begin{array}{c}\label{constantesm1}
b_2=0,\qquad c_3=0 , \qquad b_1=\,a_2,\quad   a_3=-\,c_1,\quad   b_3  =
\, c_2  ;  \\[0.3cm]
a_1 c_2=0 .
\end{array}\end{equation}
Here, the space of parameters is 3D.
\end{enumerate}
In all the  cases no more conditions are found between the structure
constants to higher orders, in
spite of  the fact that we have assumed no dependence on the deformation
parameter of the structure
constants. Note that the order by order procedure has to be solved, in
general, simultaneously for both the
quantum coproduct and the deformed commutation rules.

Now, we are in conditions to obtain the deformed commutators.
So,  the integration of the above equations to all orders gives the
general $q$--algebra of three
generators  compatible with the deformed coproduct.  We find:
\begin{enumerate}

\item
Case $\rho\neq \pm1$:
\begin{itemize}

\item[1.1)]  $c_2\neq 0$:
\begin{equation}\label{zconmutaciones01}
\begin{array}{lll}
[A,B] &=& c_1 \,\sinh (zA)/z +c_2\, B   ,\\[0.3cm]
[A,C] &=& - a_2\,\sinh (z\,\rho\, A)/z - \rho\,c_2 \,C ,\\[0.3cm]
[B,C] &=& \displaystyle \frac{a_2 \, c_1}{c_2} \frac{\sinh [z(1+\rho)
A]}{z}  \\[0.3cm]
&& \qquad+a_2\, {\rm Sym}\;[B\cosh (z\,\rho\, A)] + \rho\,c_1\, {\rm
Sym}\;[C\cosh (z A)] .
\end{array}\end{equation}

\item[1.2)]
 $c_2=0,\; a_2=0$:

\begin{equation}\label{zconmutaciones02}
\begin{array}{lll}
[A,B] &=& c_1 \,\sinh (zA)/z    ,\\[0.3cm]
[A,C] &=& 0  ,\\[0.3cm]
[B,C] &=&\displaystyle
 a_1\,\frac{\sinh [z(1+\rho)\,A]}{z(1+\rho)} + \rho\,c_1\, C  \,\cosh
(zA) .
\end{array}\end{equation}

\end{itemize}
 $(c_2=0,\; a_2=0)$.
\item
Case $\rho =+1$:
\begin{itemize}

\item[2.1)]
\begin{equation}\label{zconmutaciones1}
\begin{array}{lll}
[A,B] &=& c_1 \,\sinh (zA)/z + c_2\, B + c_3 \,C ,\\[0.3cm]
[A,C] &=& -a_2 \,\sinh (zA)/z + b_2 \,B - c_2 \,C ,\\[0.3cm]
[B,C] &=& a_1\, \sinh (2zA)/(2z) + {\rm Sym}\;\{(a_2\, B + c_1\, C)
\,\cosh (zA)\}
\end{array}\end{equation}
\end{itemize}

\item

Case $\rho =-1$:
\begin{itemize}

\item[3.1)]  $c_2\neq 0$:

\begin{equation}\label{zconmutaciones2}
\begin{array}{lll}
[A,B] &=& c_1\,\sinh (zA)/z + c_2\, B   ,\\[0.3cm]
[A,C] &=&  b_1 \,\sinh (zA)/z + c_2 \,C ,\\[0.3cm]
[B,C] &=&  {\rm Sym}\;\{(b_1\, B - c_1\, C) \,\cosh (zA)\} .
\end{array}\end{equation}

\item[3.2)]
 $c_2=0$:

\begin{equation}\label{zconmutaciones3}
\begin{array}{lll}
[A,B] &=& c_1 \,\sinh (zA)/z    ,\\[0.3cm]
[A,C] &=&  b_1 \,\sinh (zA)/z  ,\\[0.3cm]
[B,C] &=& a_1\, A +{\rm Sym}\;\{(b_1\, B - c_1\, C) \,\cosh (zA)\} .
\end{array}\end{equation}
\end{itemize}
\end{enumerate}

It is worthy to notice that in many of these case there is a form invariance of the  
commutators related with the interchange of the generators $B$ and $C$.
%%%%%%%%%%%%%%%%%%%%%%%%%%%%%%%%%%%%%%%%%%%%%%%%%%%%%%%%%%%%%%%
%%%%%%%%%%%%%%%%%%%%%%%%%%%%%%%%%%%%%%%%%%%%%%%%%%%%%%%%%%%%%%%
\sect{Equivalence and classification}\label{equivalence}

In general, two quantum algebras are said to be equivalent (isomorphic)
if there exists an
invertible (nonlinear in many cases) map between their corresponding
quantum universal enveloping
algebras as Hopf algebras. In this way, equivalence classes of quantum
deformations can be defined.
But it is clear that, due to the infinite number of possibilities given
by arbitrary nonlinear maps,
such equivalence classes are huge, and presumably some general criteria
for the choice of ``canonical"
representatives of each of them should be helpful for classification
purposes. Moreover, the
definition of a such canonical representative for the Hopf algebra
structure would  be given, indeed, in
terms of what could be properly called as ``generators" of the quantum
algebra.

In fact, the standard classification of (non deformed) Lie algebras can
be understood as a particular
application of the abovementioned procedure to their corresponding
universal enveloping algebras as
Hopf algebras, since Friedrichs theorem states that the only primitive
elements in $U(\mathfrak a)$
under the coproduct  $\co_0$ are just the generators of the Lie algebra
$\mathfrak a$. Therefore, we can define  the generators of a Lie algebra 
$\mathfrak a$ as those elements of $U(\mathfrak a)$ which have a primitive
coproduct. In this way, the (Cartan)
classification of Lie algebras is performed by obtaining appropriate
generators (primitive
elements) having the ``simplest" commutation rules (minimum number of
non-vanishing structure
constants). We stress that, in order to find such ``irreducible"
commutation rules only
equivalence transformations leaving the coproduct invariant are allowed
(in the Lie case, these are
just linear transformations).

From this perspective, we propose a definition of the ``canonical"
representatives of quantum
algebras by following a similar procedure. In this case is essential to
realize that a quantum
algebra
$U_z(\mathfrak a)$ is endowed  with a deformed coproduct
$\Delta$ which is no longer cocommutative, but $\Delta$ is (through our
quantization procedure)
invariant under
$\tilde{\sigma}$ (generalized cocommutativity).
Thus, in order to find the ``canonical" generators we shall move
within the equivalence subclass
defined through the restricted set of Hopf algebra isomorphisms of
quantum universal enveloping
algebras (always in a symmetrized basis) that leave the coproduct
formally invariant. Through such
restricted isomorphisms we shall look for representatives with
``irreducible" deformed commutation
rules having a minimal number of non-zero terms.

By proceeding in this way we have succeeded in classifying all the
non-isomorphic quantum algebras
that are contained in the three multiparameter families  given in the
previous section.
Let us explicitly
obtain them by eliminating many
irrelevant parameters through coproduct-preserving mappings.

%%%%%%%%%%%%%%%%%%%%%%%%%%%%%%
%%%%%%%%%%%%%%%%%%%%%%%%%%%%%%%%%
\subsect{Case 1:   $\rho \neq \pm 1$}

 Let us start with the case  $\rho  \neq \pm1$.
Let us consider the following transformation
\cite{dobrev}:
\begin{equation}\begin{array}{l}\label{cambio1}
\displaystyle \A=\alpha A,\quad
\B= \beta B +  \delta \frac{\sinh (A)}{z} ,\quad
\C=   \nu C +\eta \frac{\sinh (z\,\rho\,A)}{z\,\rho} ,\\[0.3cm]
\z= \alpha ^{-1} z ,\qquad
\alpha ,\; \beta ,\; \delta ,\; \nu ,\; \eta \in \mathbb C.
\end{array}\end{equation}
After this transformation the coproduct (\ref{coproductos})  becomes
 \begin{equation}\label{coproductos1}
\begin{array}{lll}
\Delta(\A) &=&  1\otimes  \A +   \A \otimes 1 , \\[0.3cm]
\Delta(\B) &=&  e^{\z \A}\otimes  \B +   \B \otimes e^{-\z \A} ,
\\[0.3cm]
\Delta(\C) &=&  e^{ \z\rho \A}\otimes  \C +   \C \otimes e^{- \z\rho
\A}  ,
\end{array}\end{equation}
i.e., it remains formally invariant. So, all the elements of the form
(\ref{cambio1}) belong of the same
class.
\smallskip

\begin{itemize}
\item[1.1)]
$c_2 \neq 0$

In this case the above mentioned change of basis (\ref{cambio1}) can be
reduced to
\begin{equation}\label{cambio11}
\begin{array}{l}\A= A / c_2,\qquad
\B=c_2 B + c_1 \frac{\sinh (zA)}{z} ,\qquad
\C=   c_2 C + a_2 \frac{\sinh (z\rho A)}{z\rho} ,\\[0.3cm]
\z= c_2 z.
\end{array}\end{equation}

\begin{itemize}
\item[1.1.1)]
Under this change of basis 
 we obtain
the new Lie commutators
\begin{equation}\label{algebra1.1.1}
[\A,\B] = \B   ,\qquad
[\A,\C] = - \,\rho\, \C ,\qquad
[\B,\C] = 0 ,
\end{equation}
which correspond to a quantization of the  Lie algebra (\ref{jacobson1}). Its
associated
bialgebra is non-coboundary, i.e., there is no classical-$r$
matrix.
\end{itemize}

\item[1.2)]
$c_2 = 0$

We have two quantum algebras. Following a procedure analogous to
(\ref{cambio11}) they can be
written as
\begin{itemize}
\item[1.2.1)]
$a_2=0, \;  c_1 \neq 0$
\begin{equation}\label{algebra1.1.2}
[\A,\B] = \sinh (\z\A)/\z    ,\qquad
[\A,\C] = 0  ,\qquad
[\B,\C] =  \rho\, \C  \,\cosh (\z\A) .
\end{equation}
The bialgebra is coboundary: the classical $r$-matrix is $r=\hat z\,\A\wedge
\B$. It is non-standard, i.e.
it verifies the classical Yang-Baxter equation. This is another quantization of the  Lie algebra
(\ref{jacobson1}).

\item[1.2.2)]
$a_2=0, \; a_1 \neq 0 ,\; c_1 = 0$
\begin{equation}\label{algebra1.2.2}
[\A,\B] = 0   ,\qquad
[\A,\C] = 0  ,\qquad
[\B,\C] =\displaystyle
  \frac{\sinh [\z (1+\rho) \A]}{ \z(1+\rho)} .
\end{equation}
This is also a coboundary deformation with standard $r$-matrix, i.e.
it verifies the modified classical Yang-Baxter equation, $r=\hat z \B\wedge\C$.\end{itemize}
\end{itemize}

Note that for $\rho=0$ we have obtained two deformations of the
extended Borel algebra and one deformation of the  Heisenberg-Weyl
algebra, both of them of Type II in
Jacobson.

%%%%%%%%%%%%%%%%%%%%%%%%%%%%%%
%%%%%%%%%%%%%%%%%%%%%%%%%%%%%%%%%
\subsect{Case 2:  $\rho =+1$}

The equivalence classes are defined by applying to
(\ref{zconmutaciones1}) the transformation:
\begin{equation} \begin{array}{l}\label{cambio2}
\A=\alpha A,\qquad
\B= \beta B + \gamma C +\delta \frac{\sinh (zA)}{z} ,\qquad
\C= \mu B + \nu C +\eta \frac{\sinh (zA)}{z} ,\\[0.3cm]
\z= \alpha ^{-1} z \qquad\qquad
\alpha ,\; \beta ,\; \gamma,\; \delta ,\; \mu ,\; \nu ,\; \eta \in
\mathbb C.
\end{array}\end{equation}

This transformation allows us to distinguish the quantum algebras
characterized by $b_2 c_3 + {c_2}^2 \neq 0$ and those  in which $b_2 c_3
+ {c_2}^2 = 0$.
\begin{itemize}

\item[2.1)] $b_2 c_3 + {c_2}^2 \neq 0$
\begin{itemize}

\item[2.1.1)]
  
\begin{equation}\label{algebra2.1.1}
[\A, \B] =\B ,\qquad [\A, \C] = -\C ,\qquad
[\B, \C ]=\frac{\sinh (2 \z \A)}{2 \z}   .
\end{equation}
This quantum algebra is just ${\cal A}_1(q)$.
The classical $r$-matrix is $r= z\,\B\wedge \C$, and is standard.

\item[2.1.2 )] 
\begin{equation}\label{algebra2.1.2}
[\A, \B] =\B ,\qquad [\A, \C] = -\C ,\qquad
[\B, \C ] = 0   .
\end{equation}
This is the complexification of
the first discovered contraction of $su_q(2)$  deformation of the
 Euclidean algebra
${\cal E}(2)$ \cite{rusos} in two dimensions. It  is a non-coboundary one.
\end{itemize}

\item[2.2)]
  $b_2 c_3 + {c_2}^2 = 0$
\medskip

\begin{itemize}

\item[2.2.1)] $a_2 c_2+b_2 c_1\neq 0$
\begin{equation}\label{algebra2.2.1}
[\A, \B] = - \frac{\sinh(\z\A)}{\z},\qquad
[\A, \C] = \B,\qquad
[\B, \C] =- {\rm Sym} \{\C  \cosh(\z \A)\}
\end{equation}
This is the symmetrized version of the well-known Jordanian
deformation of ${\cal A}(1)$ \cite{ohn} with non-standard classical $r$-matrix $r=
\z\,  \A\wedge \B$.
According
to the commutation relations (\ref{algebra2.2.1}) we obtain that
$$
{\rm Sym}(\C\,\cosh(\z\A))=\frac12(\C\,\cosh(\z\A)+\cosh(\z\A)\,\C))+
\frac{1}{12}\,\z^2\,\sinh\frac{2\z\A}{2\z},
$$
thus, the quantum algebra presented in  \cite{ohn} is a  case of
(\ref{zconmutaciones1}) with $z$-dependent parameters.

\item[2.2.2)]  $a_2 c_2+b_2 c_1=0$

This condition implies that the commutators $[\A,\B]$ and $[\A,\C]$ are proportional. 
We obtain the following 
algebras: 
\begin{itemize}
\item[2.2.2.1)]
\begin{equation}\label{algebra2.2.2}
[\A, \B] =  \frac{\sinh(\z\A)}{\z},\qquad
[\A, \C] = 0,\qquad
[\B, \C] =  \C  \cosh(\z \A)
\end{equation}
We recover a non-standard deformation of the Euclidean group  in two
dimensions $E(2)$
\cite{angel95}. The classical $r$-matrix is non-standard,
$r= \z\A\wedge \B$.

\item[2.2.2.2)]
\begin{equation}\label{algebra2.2.3}
[\A, \B] =  0,\qquad
[\A, \C] =  -\B,\qquad
[\B, \C] =    \frac{\sinh(2\z\A)}{2 \z}
\end{equation}
We have the standard deformation of $E(2)$   with  classical $r$-matrix
$r= \z \,\B\wedge \C$.

\item[2.2.2.3)]
\begin{equation}\label{algebra2.2.4}
[\A, \B] = 0
,\qquad [\A, \C] =  0
,\qquad [\B, \C] = \frac{\sinh(2\z\A)}{2 \z}
\end{equation}
It corresponds to a deformation of the Heisenberg-Weyl algebra with
classical $r$-matrix is standard, $r=\z \,\C\wedge \B$.

\item[2.2.2.4)]
\begin{equation}\label{algebra2.2.5}
[\A,\B] =  0
,\qquad [\A,\C] =  \B
,\qquad [\B, \C] = 0
\end{equation}
This is a non-coboundary deformation of the Heisenberg-Weyl algebra.
\end{itemize}
\end{itemize}
\end{itemize}

%%%%%%%%%%%%%%%%%%%%%%%%%%%%%%
%%%%%%%%%%%%%%%%%%%%%%%%%%%%%%%%%
\subsect{Case 3:  $\rho =-1$}

The classification of the quantum algebras corresponding to the case
$\rho =-1$ can be
made considering the transformation:
\begin{equation}\label{cambio3}
\A=\alpha A,\qquad
\B=\beta B + \delta \frac{\sinh (zA)}{z} ,\qquad
\C=  \nu C +\eta \frac{\sinh (zA)}{z} ,\qquad
\z= \alpha ^{-1} z ,
\end{equation}
where $\alpha ,\;\beta ,\;  \delta ,\; \nu  ,\,\eta$ are complex
numbers.

\begin{itemize}
\item[3.1)]
\ $c_2\neq 0$
\begin{equation}\label{algebra3.1.1}
[\A ,\B]=\B, \qquad [\A, \C]=\C, \qquad [ \B,\C]=0 .
\end{equation}
We have a   deformation of the dilations algebra such that there is not
$r$-matrix.

\item[3.2)]
 $c_2= 0$:
\begin{itemize}

\item[3.2.1)]
 \begin{equation}\label{algebra3.2.1}
\begin{array}{lll}
[\A,\B] &=& -  \,\sinh (\z\A)/\z    ,\\[0.3cm]
[\A,\C] &=&    \,\sinh (\z\A)/\z  ,\\[0.3cm]
[\B,\C] &=&  \, \A + (  \B +  \C) \,\cosh (\z\A) .
\end{array}\end{equation}
Like in the previous case there is not $r$-matrix.

Note that in the limit of $\z \to 0$  we recover
\[\label{algebra3.2.1limite}
[\A,\B] = -  \A    ,\qquad
[\A,\C] =   \A  ,\qquad
[\B,\C] =  \, \A + (  \B +  \C)  .
\]
that can be rewritten under a change ($ \B+\C  \to \B  $) like
\begin{equation}\label{algebra3.2.1limitec}
[\A,\B] = 0    ,\qquad
[\A,\C] =   \A  ,\qquad
[\B,\C] =  \, \A +    \B  .
\end{equation}
Therefore we have obtained a quantum deformation of the Lie algebra
(\ref{32}).

\item[3.2.2)] Other deformation of (\ref{32}), now  non-standard, is 

 \begin{equation}\label{algebra3.2.2}
[\A,\B] = 0   ,\qquad
[\A,\C] = \sinh (\z\A)/\z  ,\qquad
[\B,\C] =  \A + \B  \,\cosh (\z\A) .
\end{equation}
The $r$-matrix is   $r= \z\,\A\wedge \C$.

\item[3.2.3)]
 \begin{equation}\label{algebra3.2.3}
\begin{array}{lll}
[\A,\B] &=& - \sinh (\z\A)/\z    ,\\[0.3cm]
[\A,\C] &=&  \sinh (\z\A)/\z  ,\\[0.3cm]
[\B,\C] &=&  ( \B +  \C) \,\cosh (\z\A).
\end{array}\end{equation}
This is a deformation of the dilation algebra in two dimensions without $r$-matrix.

\item[3.2.4)]
 \begin{equation}\label{algebra3.2.4}
[\A,\B] = 0    ,\qquad
[\A,\C] =  0  ,\qquad
[\B,\C] =  \A   .
\end{equation}
It corresponds to a coboundary  deformation of the Heisenberg-Weyl
algebra with standard $r$-matrix, $r= \z\,\B\wedge \C$.

\item[3.2.5)]
\begin{equation}\label{algebra3.2.5}
[\A,\B] = -\sinh (\z\A)/\z    ,\qquad
[\A,\C] = 0  \, ,\qquad
[\B,\C] =    \C \,\cosh (\z\A) .
\end{equation}
We have a deformation of the dilatation algebra in two dimensions. In
this case the classical
$r$-matrix is non-standard,
$r= \z \A\,\wedge \B$.
\end{itemize}
\end{itemize}

%%%%%%%%%%%%%%%%%%%%%%%%%%%%% SECTION 4 %%%%%%%%%%%%%%%%%%%%%%%%%%%%%%%
\sect{Conclusions and remarks}\label{conclusiones}
%%%%%%%%%%%%%%%%%%%%%%%%%%%%%%%%%

This quantization method
can be simultaneously applied
and successfully solved for a multiparameter family of Lie bialgebras that share some
structural properties.

Throughout the paper we have considered certain 3D complex Lie
bialgebras. In particular,
the parameter $\rho$ is complex but we do  also have   isolated
solutions for $\rho=\pm 1$.
On the other hand, the comparison with the complete classification of 3D
real Lie bialgebras given in
\cite{gomez} can be worked out by considering the isomorphisms among the
complexified
versions of 3D real Lie algebras.

In this way it can be shown that all the quantum algebras that we have
obtained in section 4 are
quantizations of the complexifications of the dual version of the Lie
bialgebras given in \cite{gomez}.
In order to find out the correspondence explicitly, in
Table III of
\cite{gomez} we have to identify a given algebra $\mathfrak a$ with the
complex version of the
{\it dual} Lie algebra $\mathfrak g^\ast$ and, consequently, the dual of
the cocommutator $\eta$ will
have to be isomorphic to one of the algebras $\mathfrak g$ in the first
row of such Table. By
proceeding in this way we find the following correspondences (we write
first the Lie bialgebras
$(\mathfrak g^\ast,\mathfrak g\equiv
\eta^\ast)$ as labeled in Table III of
\cite{gomez} and afterwards the corresponding quantum algebra according
to our classification):
\begin{itemize}

\item[ ]
5   $\; \rightarrow \;$  1.2.2) ;
6   $\; \rightarrow \;$  1.2.1) ;
7   $\; \rightarrow \;$  1.1.1) ;

\item[ ]
(1)   $\; \rightarrow \;$  2.1.1) ;
(2), (4)   $\; \rightarrow \;$  2.1.1) ;
(3)   $\; \rightarrow \;$  2.2.1) ;
9   $\; \rightarrow \;$  2.1.2) ;
11, 11'   $\; \rightarrow \;$  2.2.2.2) ;

\item[ ]
10  $\; \rightarrow \;$  2.2.2.4) ;
5$_{\rho=1}$ $\; \rightarrow \;$  2.2.2.3) ;
6$_{\rho=1}$  $\; \rightarrow \;$  2.2.2.1)$ $ ;
7$_{\rho=1}$  $\; \rightarrow \;$  2.1.2)$_{\rho=1}$ ;

\item[ ]
 5'  $\; \rightarrow \;$  3.2.4) ;
 8  $\; \rightarrow \;$  3.2.1) ;
(14)  $\; \rightarrow \;$  3.2.2) ;
(11)  $\; \rightarrow \;$  3.2.3) ;
6$_{\rho=-1}$  $\; \rightarrow \;$  3.2.5)$_{\rho=-1}$ ;
7$_{\rho=-1}$  $\; \rightarrow \;$  3.1)$_{\rho=-1}$ .

\end{itemize}

In this way we can realize that our choice (\ref{coconmutador}) for the
cocommutator implies that we
have just obtained the quantizations for the full  set of dual Lie
bialgebras
$(\mathfrak g^\ast,\mathfrak g\equiv
\eta)$ of
\cite{gomez} such that $\eta^\ast\equiv{\mathfrak r}_3(\rho)$ for all
values of $\rho$. In fact, the
$\rho$ parameter in (\ref{coconmutador}) is identified with the one
appearing in Gomez's
classification.

%%%%%%%%%%%%%%%%%%%%%%%% ACKNOWWLEDGMENTS %%%%%%%%%%%%%%%%%%%%%%%%%%%%%
\section*{Acknowledgments}
%%%%%%%%%%%%%%\%%
This work has been partially supported by
DGI of the  Ministerio de Ciencia y Tecnolog\'{\i}a
(Projects BMF2002-0200 and BFM2000-1055), the FEDER Programme  and
Junta de Castilla y Le\'on (Projects VA085/02 and BU04/03) (Spain).
The  visit of E.C. to Valladolid
have been  financed by  Universidad de Valladolid, by CICYT-INFN and by
Ministerio
de Educaci\'on y Cultura (Spain).

%%%%%%%%%%%%%%%%%%%%%% BIBLIOGRAPHY %%%%%%%%%%%%%%%%%%%%%%%%%%%

%%%%%%%%%%%%%%%%%%%%
\end{document}